\documentclass[12pt,a4paper,oneside]{amsart}
\usepackage{amsfonts, amsmath, amssymb, amsthm, amscd, hyperref}
\usepackage[T2A]{fontenc}
\usepackage[cp1251]{inputenc}
\usepackage[english]{babel}

\usepackage{cmap}
\usepackage[english]{babel}
\usepackage{amssymb}
\usepackage{amsmath}
\newtheorem{lemma}{Lemma}

\newtheorem{fact}{Theorem}
\newtheorem{corollary}{Corollary}
\newtheorem{question}{Question}
\newtheorem{example}{Example}

\theoremstyle{definition}
\newtheorem{definition}{Definition}
\newtheorem{remark}{Remark}
\newcommand{\la}{\lambda}

\newcommand{\e}{\varepsilon}
\newcommand{\al}{\alpha}

\newcommand{\men}{\leqslant}
\newcommand{\bol}{\geqslant}
\newcommand{\bra}{\langle}
\newcommand{\ket}{\rangle}
\newcommand{\B}{B}
\newcommand{\R}{\mathbb{R}}

\newcommand{\norm}[1]{\left\| #1 \right\|}
\newcommand{\ovl}[1]{\overline{#1}}
\def\vn{\mathop{\rm int}}
\def\co{\mathop{\rm co}}
\def\Lin{\mathop{\rm Lin}}
\def\diam{\mathop{\rm diam}}
\def\dim{\mathop{\rm dim}}

\begin{document}

\title{Convex hull deviation and contractibility}

\author{G.M. Ivanov}

\address{Department of Higher Mathematics, Moscow Institute of Physics and Technology,  Institutskii pereulok 9, Dolgoprudny, Moscow
region, 141700, Russia}
\address{
National Research University Higher School of Economics,
School of Applied Mathematics and Information Science,
Bolshoi Trekhsvyatitelskiy~3, Moscow, 109028, Russia}
 \email{grimivanov@gmail.com}

\maketitle
\begin{abstract}
We study the Hausdorff distance between  a set   and its convex hull.
Let $X$ be a Banach space, define  the CHD-constant of space $X$ as the supremum of this distance for all subset of the unit ball in $X$.
In the case of finite dimensional Banach spaces we obtain the exact upper bound of the CHD-constant depending on the dimension of the space.
We give an upper bound for the CHD-constant in $L_p$ spaces.
We prove that CHD-constant is not greater than the maximum of the Lipschitz constants of metric projection operator onto  hyperplanes.
This implies that for a Hilbert space  CHD-constant equals 1. 
We prove criterion of the Hilbert space  and study the contractibility of  proximally smooth sets in uniformly  convex and uniformly smooth Banach spaces.
\end{abstract}
%\begin{keyword}
%convex hull deviation, Hilbert space criterion, contractibility, proximally smooth sets
%\end{keyword}
%\end{frontmatter}

\section{Introduction}
Let $X$ be a Banach space. For a set $A \subset X$ by $\partial A, \vn A$ and $ \co A$
we denote the boundary, interior and convex hull of $A,$ respectively.
We use $\bra p,x \ket$ to denote the value of functional $p \in X^*$ at the vector $x \in X.$
For $R>0$ and $c \in X$ we denote by $\B_R(c)$ a closed ball with center $c$ and radius $R.$

By  $\rho(x, A)$ we  denote  distance between the point  $x\in X$ and set $A.$ 
We define the  deviation from set $A$ to set $B$ as follows
\begin{equation}
	h^{+}(A,B) = \sup\limits_{x \in A} {\rho(x, B)}.
\end{equation}
In case $ B \subset A $, which
takes place below, the deviation $	h^{+}(A,B)$ coincides with the Hausdorff  distance
 between the %todo здесь нужен предлог? 
 sets $ A $ and $ B $.

Given $D\subset X$ the deviation $h^+(\co D,D)$ is called the {\it convex hull deviation} (CHD) of $D$.

We define {\it CHD-constant} $\zeta_X$ of $X$ as
$$\zeta_X = \sup_{D \subset \B_1(o)} h^+(\co D,D).$$
\begin{remark} \label{rem_1}
 Directly from our definition it follows that for any normed linear space $X$ 
we have $1 \men \zeta_X \men 2.$
\end{remark}
We denote by  $\ell_{p}^{n}$ the  $n$-dimensional real vector space with the $p$-norm.

This article contains estimates  for the CHD-constant for different spaces and some of its geometrical applications. 
In particular, for finite-dimensional spaces we obtain the exact upper bound of the CHD-constant depending on the dimension of the space:
\begin{fact}\label{UVOfinite}
Let $X_n$ be  a normed linear space, $\dim X_n = n \bol 2,$ then $\zeta_{X_n} \le
2\frac{n-1}{n}.$ 
If 
$X_n=\ell_{1}^{n}$ or $X_n=\ell_{\infty}^{n},$
then the estimate is reached.
\end{fact}

Let the sets $P$ and $Q$  be  the intersections of the unit ball with two parallel affine hyperplanes of dimension $k$ and $P$ is a central section.
In Corollary \ref{cor1} we obtain the exact upper bound of the homothety coefficient, that provides covering of $Q$ by $P.$

The next theorem gives an estimate for the CHD-constant in the $L_p, 1 \men p \men +\infty $ spaces:
\begin{fact}\label{th_Lp}
For any $p \in [1,\, +\infty]$
\begin{equation}\label{L_p_modul}
\zeta_{L_p} \men 2^{\left| \frac{1}{p} - \frac{1}{p{'}} \right|},
\end{equation}
where $\frac{1}{p} + \frac{1}{p{'}}=1.$
\end{fact}

Theorem \ref{th_chi} shows that CHD-constant is not greater than the maximum of the Lipschitz constants of metric projection operator onto  hyperplanes.
This implies that for Hilbert space  CHD-constant equals 1.
Besides that, we prove the criterion of a Hilbert space in terms of CHD-constant. 
The  idea of the proof is analogous to the idea used by A. L. Garkavi in  \cite{Garkavi}.
\begin{fact}\label{UVOcriterion}
The equation   $\zeta_X =1$ holds for a Banach space $X$ iff $X$ is an Euclidian space or 
$\dim X= 2$.
\end{fact}

	In addition we study the contractibility of a covering of the convex set with balls. 
\begin{definition}
	A covering of a convex set with balls is called {\it admissible} if it consists of a finite number of balls with centers in this set and the same radii.
\end{definition}
\begin{definition}
A family of balls is called {\it admissible} when it is an admissible covering of the convex hull of its centers.
\end{definition}
We say that a covering of a set by balls is conractible when the union of these balls is contactible.
It is easy to show that in two-dimensional and Hilbert spaces any admissible covering  is contractible 
(see Lemmas \ref{lemma_collapse} and \ref{lemma_collapse2}). On the other hand, using Theorem \ref{UVOcriterion}, we  prove the following statement.
\begin{fact} \label{ballcriterion}
	In a three dimensional Banach space $X$ every admissible covering is contractible iff $X$ is a Hilbert space. 
\end{fact}
For  3-dimensional spaces we consider an example of  an admissible covering of a convex set with four balls that is not contractible.
To demonstrate the usefulness of this technics in Theorem \ref{th_prox} we obtain the sufficient condition
 for the contractibility of the proximally smooth sets 
in uniformly convex and uniformly smooth Banach space.  
%%%%%%%%%%%%%%%%%%%%%%%%%%%%%%%%%%%%%%%%%%%%%%%%%%%%%%%%%%%%%%%%%%%%%%%%%%%%%%%%%%%
\section{Proof of Theorem 1 and some other results}
\begin{lemma}\label{lemma o svyaznosti}
  Suppose the set $\B_1(o)\setminus \vn \B_r(o_1)$ is nonempty. Then it is arcwise connected.
\end{lemma}
{\bf Proof.}\\
We suppose that $o \neq o_1,$ otherwise the statement is trivial. 
Let $z$ be the point of intersection of  ray $o_1o$ and the boundary of the closed ball $\B_1(o)$. 
The triangle inequality tells us that $\B_1(o)\setminus \vn \B_r(o_1)$ contains $z$ (because the set is nonempty). 
We claim that $\partial \B_1(o) \setminus \vn \B_r(o_1)$ is arcwise connected and thus prove the lemma. 
It suffices to show that in the two dimensional case every point of $\partial \B_1(o) \setminus \vn \B_r(o_1)$ 
is connected with $z.$
 Suppose, by contradiction, that it is not true. 
This means that 
there exist points $a_1, b_1 \in \partial \B_r(o_1) \cap \partial \B_1(o)$ lying on the same side of the line $oo_1$
such that the arc  $a_1b_1$  of the circle $\partial \B_1(o)$ contains a point $c \notin \B_r(o_1),$ that is $\norm{c- o_1} > r.$

Consider two additional rays $oa$ and $ob$ codirectional with $o_1a_1$ and $o_1b_1$ respectively, where $a, b \in \partial \B_1(o). $ 
Since balls $\B_1(o)$ and $\B_r(o_1)$ are similar,  we have $a_1b_1 \parallel ab.$ So, the facts that points $a,b,a_1,b_1$ lie on the same side of $oo_1$ line, $ oa \cap o_1a_1 = \emptyset, ob \cap o_1b_1 = \emptyset$ and that a unit ball is convex, imply that segments $ab$ and $a_1b_1$ lie on the same line, this contradicts $\norm{c- o_1} > r.$

{\bf Proof of  Theorem \ref{UVOfinite}.}\\
Denote $r_n = 2 \frac{n-1}{n}$.

Suppose the inequality doesn't hold. It means that there exists a Banach space $X_n$ with dimension $n \bol 2$, a set $D \subset \B_1(o) \subset X_n$ and a point $o_1 \in \co D$, such that $\B_{r_n}(o_1) \cap D = \emptyset$. But if $o_1 \in \co D $, then $o_1 \in \co (\B_1(o) \backslash \vn \B_{r_n}(o_1))$. 
According to Lemma \ref{lemma o svyaznosti} the set $B = \B_1(o) \backslash \vn \B_{r_n }(o_1)$ is connected.
 So, taking into consideration the generalized Caratheodory's theorem (\cite{RockWetsVarAn}, theorem 2.29), we see that the point $o_1$ is a convex combination of not more than $n$ points from $B$. 
These points denoted as $a_1, \cdots, a_k, \; k \men n$, may be regarded as vertices of a $(k -1)$-dimensional simplex $A$ and point $o_1 = \al_1 a_1 + \cdots + \al_k a_k$ lies in its relative interior 
($\al_i > 0, \al_1 + \cdots + \al_k = 1$).

Let $c_l$ be the point of intersection of ray $a_l o_1$ with the opposite facet of the simplex 
$A.$ 
So, $o_1 = \al_l a_l + (1 - \al_l)c_l.$ 
Then 
$$\|o_1-a_l\|= (1 - \al_l)\|c_l-a_l\|.$$ 
And $[c_l,a_l] \subset A \subset \B_1(o)$ implies that $\|a_l - c_l\| \men 2,$ 
for all $l \in \ovl{1,k}.$ 
Therefore $r_n < \|o_1 - a_l\| \men 2(1 - \al_l).$ Thus $\al_l < 1 - \frac{r_n}{2} < \frac{1}{n},$ and finally
$\al_1 + \cdots + \al_k < \frac{k}{n} \men 1.$ Contradiction.

Now let us show that the estimate is attained for spaces $\ell_1^n, \, \ell_{\infty}^n$.

Consider $\ell_{1}^n.$ Let $A = \{e_i\}^n_{i = 1}$ be a standard basis for $\ell_1^n$ space and $b =\frac{1}{n}(e_1 + \ldots + e_n) \in \co \{e_1, \ldots, e_n\}$. The distance between point $b$ and an arbitrary point from $A$ is $\|a_i-b\| = 2 \frac{n-1}{n}$.

Consider $\ell_{\infty}^n$. 
Let $a_{ij}=(-1)^{\delta_{ij}},$ where $\delta_{ij}$ is Kroneker symbol, $a_i= (a_{i1}, \cdots, a_{in})$ and $A = \{a_i\}^n_{i = 1}.$ 
Now let $b  = \frac{1}{n}(a_1 + \ldots + a_n)= \left(\frac{n-2}{n}, \cdots, \frac{n-2}{n} \right) \in \co \{a_1, \ldots, a_n\}.$ 
And the distance from point $b$ to an arbitrary point from $A$ is $\|a_i-b\| = 2 \frac{n-1}{n}$.

$\blacksquare$

So, Theorem 1 and inequality $\zeta_X\ge 1$ imply the CHD-constant of any 2-dimensional normed space equals $1$. Obviously, CHD-constant of $\ell_{1}$ space equals~$2$.

\begin{remark}\label{remark}
Let $X$ be a Banach space, $\dim X = n.$ Then for every $d < \zeta_X$ there exists a set $A$ that consists of not more than $n$ points and meets the condition $h^+(\co{A},A)=d.$
\end{remark}

%%%%%%%%%%%%%%%%%%%%%%%%%%%%%%%%%%%%%%%%%%%%%%%%%%%%%%%%%%%%%%%%%%%%%%%%%%%%%%%%%%%%%%%%%%

\begin{corollary}\label{cor1}
  Let sets $P$ and $Q$ be plane sections of   the unit ball with two parallel hyperplanes of dimension $k,$ 
	and let the hyperplane containing $P$ contains $0$ as well. 
	Then it is possible to cover $Q$ with the set $\min\{2\frac{k}{k+1}; \zeta_X\}P$ using parallel translation.
\end{corollary}
{\bf Proof.}\\ \noindent
Define $\eta = \min\{2\frac{k}{k+1}; \zeta_X\}.$
	Due Helly theorem it sufficies to prove that we could cover any $k$-simplex $\Delta \subset Q$
	with the set $\eta P.$
	
	Let us consider $k$-simplex $\Delta \subset Q$ with verticex $\{x_1, \cdots, x_{k+1}\}.$
	Due to the definition of the $\zeta_X$ and by Theorem \ref{UVOfinite} for any set of indices
	$I \subset \overline{1, (k+1)},$ 
	we have  $\co\limits_{i \in I}\{x_i\} \subset \bigcup\limits_{i \in I} (\B_\eta (x_i) \cap \Delta).$
	Using KKM theorem \cite{KKM_original} we obtain that  
	$S = \bigcap\limits_{i \in \overline{1, (k+1)}} (\B_\eta (x_i) \cap \Delta) \neq \emptyset.$
 	Then $ \Delta  \subset \B_\eta(s)$, where $s \in S \subset \Delta.$
$\blacksquare$\\

Let us show that Hilbert and 2-dimensional Banach spaces meet the requirements of Theorem \ref{ballcriterion}.
We consider the area covered with balls to be shaded. Balls' radii may be taken equal to 1. 

\begin{lemma} \label{lemma_collapse}
	Let $X$  be a Banach space, $\dim X = 2,$ then any admissible covering is contractible.
\end{lemma}
{\bf Proof.}\\
Without loss of generality, let we have an admissible covering of a convex set $V$ by balls $\B_1(a_i), \, i=\ovl{1,n}.$
Let us put $S = \bigcup\limits_{i\in \ovl{1,n}}\B_{1}(a_i).$
Since the unit ball is a convex closed body,  the set
$S$ is homotopy equivalent to  its nerve \cite{Nerv}, in our case it is  finite CW complex.
Therefore, $S$ is contractible iff $S$ is connected, simply connected and its homology groups $H_k(S)$ are trivial for $k \bol 2.$
Obviously, $S$ is connected set.

Let us show that the set $S$ is simply connected and $H_k(S) = 0$ for $k \bol 2$.
The unit circle is a continuous closed line without self-intersections, 
it divides a plane in two parts.
 A finite set of circles divides a plane in a finite number of connected components. Let us now shade the unit balls. 

 It is remarkable, that the problem is stable against subtle perturbations of norm. To be more precise: if a norm does not meet the requirements of the theorem, then  there exists a polygon norm, which does not meet them too.

%Consider  the norm $\norm{\cdot}_1$ on a plane with a ball $\B_1(o)$ and an admissible covering of a convex set $V$, 
%which consists of the balls $\B_1(a_i), \, i=\ovl{1,n}$, that is not contractible. 
Let us choose a bounded not-covered area $U$ with shaded boundary. 
It is possible to put a ball of radius $3\e_1 \,(\e_1 >0)$ inside this area.
There exists $\e_2 \,(\e_2 >0)$ such that if $\B_1(a_{i_1})\cap \B_1(a_{i_w})  = \emptyset$ for  $i_1, i_2 \in \ovl{1,n},$ 
than   $\B_{1+\e_2}(a_{i_1})\cap \B_{1+\e_2}(a_{i_2})  = \emptyset.$
Denote $\e = \min\{\e_1, \; \e_2\}.$

Consider the following set
 $$\B_{1}^c(o)= \bigcap_{p \in C}\{x: \bra p, x \ket \men 1\},$$ 
where $C$ is a finite set of unit vectors from space $X^*$, such that $C = - C$. 
So, $\B_1^c(o)$ is the unit ball for some norm.
According to \cite{PolBal}, Corollary 2.6.1, it is possible to pick such a set $C$, 
that $h^{+}\left(\B_1^c(o), \B_1(o) \right) \men \e.$ Then the set of balls $\B_{1}^c(a_i), \, i=\ovl{1,n}$ is admissible covering, contains the boundary of $U$, because $\B_1(o) \subset \B_{1}^c(o),$ and it does not cover $U$ entirely.
Furthermore nerve, and consequently homology group, of the sets $\bigcup\limits_{i\in \ovl{1,n}}\B_{1}^c(a_i)$ and $S$ are coincide.

Now it suffices to  show that the statement of the lemma is true in case of a polygon norm. 
In this case $S$ is the neighborhood retract in $\R^2$ (see \cite{SandersonRorkePiceLinTop}), 
therefore straightforward from Alexander duality (see \cite{Dold}, Chapter 4, \S 6) we obtain that $H_k(S) = 0$ for $k \bol 2.$

Now we shall prove that $S$ is simply connected.
Assume  the contrary, there exist a norm, an admissible covering of a convex set $V$ by balls $\B_1(a_i), \, i=\ovl{1,n}$ 
and non-shaded bounded set $U$ with a shaded boundary. 
Note that its boundary appears to be a closed polygonal line without self-intersections. 
Let  us define set $A=\co \{ a_i | i=\ovl{1,n}\}.$

Let $x$ be an arbitrary point of the set $U$. 
The  union of the balls $\B_1(a_i)$ is admissible covering of the set $A$, thus $x\not\in A.$ Then there exists a line $l_a$ that separates $x$ from set $A$. This line may serve as a supporting line of set $A$. 
Let $l\parallel l_a$ be a supporting line of $U$ in a point $v$, such that sets $U$ and $A$ lie at one side from line $\ell.$
 Line $l$ divides the plane in two semiplanes.
Let $H_+$ be the semiplane that does not contain $A$, 
 we denote the other semiplane as $H_-.$ Let points $p,q \in l$ lie on different sides from $v.$ We want to choose all the edges of polygonal curve $\partial U$, that contain point $v.$ We will call them $vb_i, i \in \ovl{1,k}:\; \cos{\angle pvb_i} > \cos{\angle pvb_j},\; i>j.$

Note that it is impossible for any of the edges to lie on line $l.$ 
Otherwise $l$ is supporting line for a ball $\B_1(a_p), p \in \ovl{1,n},$ 
and $\B_1(a_p) \cap H_+ \neq \emptyset,$ so we come to the contradiction.
 We may pick such a number $\e$ that the ball $\B_{\e}(v)$
 intersects only with particular edges of polygonal curve $\partial U$. 
From now on we use $p,\, q, \, b_i, i \in \ovl{1,k}$ for points of intersection of circle $\partial \B_{\e}(v)$ with corresponding edges.
Since $v \in \partial U,$ it follows that there  exists a point $z$ on circle $\partial \B_{\e}(v)$,
 such that the interior of segment $vz$ lies in $U$ 
and the ray $vz$ lies between $vb_1$ and $vb_k.$ 
Then, since the ball is convex, there is no such  ball $\B_1(a_i)$, 
that simultaneously covers a point from the interior of $vb_1$ 
and a point from $vb_k,$ i.e. point $v$ is covered by at least two balls,
 and the centers of these balls  $a_i, \, a_j$ are divided by ray $vz$ in semiplane $H_-.$
Again, since the ball is convex, point $x = vz \cap a_ia_j$ is not covered by balls  $\B_1(a_i), \; \B_1(a_j),$ thus 
$\|a_i-a_j\|=\|x-a_i\|+\|x-a_j\|>2,$ which contradicts  the fact that $a_i$ and $a_j$ are contained in ball $\B_1(v)$.
	$\blacksquare$
%%%%%%%%%%%%%%%%%%%%%%%%%%%%%%%%%%%%%%%%%%%%%%%%%%%%%%%%%%%%%%%%%%%%%%%%%%%%%%%%%%%%%%%%%%%%%%%%%%%%%%%%%%
\begin{lemma} \label{lemma_collapse2}
	Let $X$ be an Euclidean space. Then any admissible covering is contractible.
\end{lemma}
{\bf Proof.}\\
Let us remind that a closed convex set is contractible and  in a Hilbert space the projection onto a closed convex set is unique. 
Since a projection onto a convex set is a continuous function of the projected point,
 it is enough to prove that a line segment, 
which connects a shaded point with its projection onto a convex hull of centers of an admissible covering,
 is shaded. 
Suppose that we have an admissible set of balls.
 The convex hull of its center is a polygon, 
Let us call it $C.$ 
If a shaded point $a$ is projected onto the $v$-vertex of the polygon, 
then  the segment $av$ is shaded as well.
 Let a shaded point $a$ lying in the ball $\B_1(v)$ from a set of balls be projected onto the point $b \neq v.$
Let $L$ be a hyperplane passing through point $b$ and perpendicular to the line segment $[a,b].$ It divides the space in two half-spaces. The one with the point $a$ we call $H_a$. $C$ is convex, thus it contains the segment $[v,b]$. Then it is impossible for point $v$ to lie in $H_A$, so $\angle{abv} \bol \frac{\pi}{2},$ i.e. $\norm{v-a} \bol \norm{v-b}.$ Thus, $b \in \B_1(v)$ and, consequently,  $ ab \subset \B_1(v)$. 
	$\blacksquare$\\

%%%%%%%%%%%%%%%%%%%%%%%%%%%%%%%%%%%%%%%%%%%%%%%%%%%%%%%%%%%%%%%%%%%%%%%%%%%
\section{Upper bound for CHD-constant in a Banach space}

Let $J_1(x)=\{p \in X^*\mid \bra p, x \ket = \norm{p} \cdot \norm{x}=\norm{x}\}.$
Let us introduce the following characteristic of a space:
$$
    \xi_X = \sup\limits_{{\norm{x}=1,} \atop {\norm{y}=1}}\sup\limits_{p \in J_1(y)}{\norm{x - \bra p, x \ket y}},
$$

Note that if $y \in \partial \B_1(0), \;p \in J_1(y),$ 
then vector $(x - \bra p,x\ket y)$ is a metric projection of $x$ onto the hyperplane $H_p = \{x \in X: \bra p, x \ket = 0\}.$ 
So, $\xi_X = \sup_{y \in \B_1(o)}\sup_{p \in J_1(x)}{\xi_X^p},$ 
where $\xi_X^p$ is  half of diameter of a unit ball's projection onto the hyperplane $H_p.$ 
This implies the following remark.

%\begin{remark} \label{rem_2}
%$\xi_X$  is the maximum of the  Lipschitz constants of the metric projection operator onto  hyperplanes.
%\end{remark}

Let us use $\xi_X$ for estimation of CHD-constant of $X$:
\begin{lemma}\label{lemma o schenie}
  Let $y \in \co{\left[\B_1(o) \backslash {\vn{\B_r(y_1)}}\right]}$ and let $p \in J_1(y).$
	There is $x \in \B_1(o)\backslash \vn{\B_r(y_1)}$ such that $\bra p,x \ket = \bra p, y \ket.$
	
	Then in hyperplane $H_p = \{x \in X:  \bra p, x \ket= \bra p, o_1 \ket\}$ there exists a point $x$, such that $x \in \B_1(o) \backslash \vn{\B_r(o_1)}.$
\end{lemma}
{\bf Proof.}\\
Define set $B = \B_1(o) \backslash {\vn \B_r(y)}.$ Since $y \in \co{B},$ there exist points $a_1, \cdots, a_n \in B$ 
and a set of positive coefficients $\la_1, \ldots, \la_n \;  (\la_1 + \ldots + \la_n = 1)$, such that 
    \begin{equation} \label{UVO lemma 2}
        y = \la_1 a_1 + \ldots +\la_n a_n.
    \end{equation}
Let $H_p^+ = \{x \in X: \bra p, x \ket \bol \bra p, y \ket.$ 
According to Lemma \ref{lemma o svyaznosti} set $B$ is connected, thus, since $B \backslash H_ p^+$ is not empty, 
if the statement we prove is not true, 
we arrive at $B \cap H_p^+ = \emptyset.$ 
Then $\bra p, a_i \ket < \bra p, y \ket$ and formula (\ref{UVO lemma 2}) implies 
$$
	\bra p, y \ket = \la_1 \bra p, a_1 \ket + \ldots + \la_n \bra p, a_n \ket < \bra p, y \ket.
$$
Contradiction.$\blacksquare$\\
\begin{lemma} \label{zeta cherez chi}
  \begin{equation}
  \zeta_X \men \sup\limits_{\norm{y}=1} \inf\limits_{p \in J_1(y)} \sup\limits_{x \in \B_1(o): \atop \bra p, x - y\ket = 0} {\|x - y\|}
 \end{equation}
\end{lemma}
{\bf Proof.}\\
%Второе неравенство очевидно, так как $$\xi_X = \sup_{y \in \B_1(O), p
%\in J_1(y)}\sup_{x \in \B_1(O): \bra p, x - y\ket = 0} {\|x-y\|}.$$
Let $\varepsilon$ be a positive real number.
 Then, according to the definiton of the CHD-constant, 
there exists set $D \subset \B_1(o)$, 
such that $h^+(\co D, D) \bol \zeta_X - \varepsilon.$ 
It means that there exists point $y \in \co D:$ $\rho(y, D) \bol \zeta_X - 2 \varepsilon.$ 
Let us put $r = \rho(y, D).$ So, $D \subset \B_1(o)\setminus \vn \B_r(y).$ 
Hence, $y \in \co{[\B_1(o)\setminus \vn \B_r(y)]}.$ 
Now  let $p \in J_1(y).$

According to Lemma \ref{lemma o schenie} 
there exists vector $x \in \B_1(o)\setminus \vn \B_r(y): \bra p, x - y\ket = 0.$ And $r \men \|x - y\|.$ 
Therefore, $\zeta_X \men \rho(y,\, D) + 2\e = r+ 2\e \men \|x-y\| + 2 \varepsilon.$ Now let $\varepsilon$ tend to zero. The lemma is proved. 
$\blacksquare$\\
It becomes obvious that 
$$
\xi_X = \sup\limits_{y \in \B_1(o), \atop p
\in J_1(y)}\sup_{x \in \B_1(o): \atop
\bra p, x - y\ket = 0} {\|x-y\|}.
$$
Then Lemma \ref{zeta cherez chi} implies
\begin{fact}\label{th_chi}
  $\zeta_X \men \xi_X.$
\end{fact}
Using Remark \ref{rem_1} and Theorem \ref{th_chi} we get
\begin{corollary} \label{Hilbert}
  If $H$ is a Hilbert space, then $\zeta_H = 1.$
\end{corollary}
With the following lemma we can pass to finite subspace limit in CHD-constant calculations.
\begin{lemma}\label{lemma_o_zamiuk}
Let $X$ be a Banach space and $\{x_1, \; x_2, \; \cdots\}$ be a vector system in it, such that the subspace $\check{X}=\Lin{\{x_1, \; x_2, \; \cdots\}}$ is dense in $X$.
Then
  \begin{equation}\label{infinit_lin1}
        \zeta_X= \lim_{n \to \infty}{\zeta_{X_n}},
  \end{equation}
	where $X_n = \Lin\{x_1, \cdots, x_n\}.$
\end{lemma}
{\bf Proof.}\\
Let us set $ \zeta =\zeta_X,$ and fix  a real number $\e > 0.$
 Since $X_n \subset X_{n+1}\subset X $, the sequence $\zeta_{X_n}$ is monotone and bounded and, consequently, convergent. 
Let $\zeta_2 = \lim\limits_{n \to \infty}{\zeta_{X_n}}.$ 
Since $X_{n}\subset X$ 
it follows that $\zeta_2 \men \zeta.$ 
According to the CHD-constant definition there exists a set $A \subset \B_1(o)$ 
and a point $d \in \co{A}$, such that $\rho(d,A) > \zeta - \frac{\e}{2}.$ 
Since $d \in \co{A},$ there exist a natural number $N,$
points $a_i \in A,$ and numbers  $\al_i \bol 0,\;  i \in \ovl{1,N}, \, \al_1 + \cdots + \al_N=1,$ such that $d = \al_1 a_1 + \cdots + \al_N a_N. \; $

Then $\norm{d - a_i}> \zeta - \frac{\e}{2}, \; i \in \ovl{1,N}$ 
Since $\ovl{\check{X}} = X,$ it is possible to pick points $b_i\in \B_1(o) \cap \check X, \; i \in \ovl{1,N}$, 
so that $\norm{a_i - b_i}\men \frac{\e}{4}.$ 
According to the definition of a linear span for some natural $n_i$ we have:
 $b_i \in X_{n_i}.$ Let $M = \max{n_i}, i \in \ovl{1,N}.$ 
Consider set $B = \{b_1, \; \cdots, \; b_N\}$ in the space $X_M$.
 Let $d_{\e} = \al_1 b_1 + \cdots + \al_N b_N \in \co{B},$ 
then 
$$\norm{d_{\e} - d}= \norm{\sum_{j=1}^{N}{\al_j (b_j - a_j)}} \men \sum_{j=1}^{N}\al_j \norm{b_j - a_j}\men \frac{\e}{4},$$ 
so for every $i\in\ovl{1,N}$ we have
$$
\norm{d_{\e} - b_i}=\norm{(d_{\e}-d) +(d - a_i) + (a_i -b_i)} \bol \norm{d - a_i} - \norm{d_{\e}-d} - \norm{a_i -b_i} \bol \zeta - \e.
$$

Thus $\zeta - \e \men h^+(\co{B},B)\men \zeta_{X_M} \men \zeta_2 \men \zeta,$ and since $\e>0$ was chosen arbitrarily, $\zeta= \zeta_2.$
$\blacksquare$
%%%%%%%%%%%%%%%%%%%%%%%%%%%%%%%%%%%%%%%%%%%%%%%%%%%%%%%%%%%%%%%%%%%%%%%%%%%%%%%%%%%%

Let $p' \in [1;\; +\infty]$ be such that $\frac{1}{p}+\frac{1}{p'}=1,$
$r = \min\{p,p'\}, \; r'=\max\{p, p'\}.$
\begin{lemma}
	Given $p \in [1,\, +\infty].$ Let $x_i \subset L_p, 1 \men p \men \infty, i=1, \cdots, k;$
	$$\sum_{i=1}^{k} {\al_i}=1, \; \al_i \bol 0 \;(i=1, \cdots, k) , \;x_0 = \sum_{i=1}^{k} {\al_i x_i}.$$
	Then
\begin{equation} \label{inequality1}
		\left(\sum_{i=1}^{k}{\al_i \norm{x_i - x_0}_p^r}\right)^{\frac{1}{r}} \men
    2^{-\frac{1}{r'}} \left( \sum_{i=1,j=1}^{k}{\al_i \al_j \norm{x_i -x_j}_p^r}\right)^{\frac{1}{r}},
	\end{equation}
    \begin{equation} \label{inequality2}
    \left( \sum_{i=1,j=1}^{k}{\al_i \al_j \norm{x_i -x_j}_p^r}\right)^{\frac{1}{r}} \men
    2^{\frac{1}{r}}\max_{1 \men i \men k}{\norm{x_i}_p}.
	\end{equation}
	
If $1 \men p \men 2,$ then the latter inequality can be strengthened:
    \begin{equation} \label{inequality3}
    \left( \sum_{i=1,j=1}^{k}{\al_i \al_j \norm{x_i -x_j}_p^r}\right)^{\frac{1}{r}} \men
    2^{\frac{1}{r}} \left( \frac{k-1}{k} \right)^{\frac{2}{p}-1} \max_{1 \men i \men k}{\norm{x_i}_p}.
	\end{equation}
\end{lemma}
{\bf Proof.}\\
The inequality (\ref{inequality2}) follow from Schoenberg's inequalities (\cite{Williams}, Theorem 15.1):
$$
\left( \sum_{i=1,j=1}^{k}{\al_i \al_j \norm{x_i -x_j}_p^r}\right)^{\frac{1}{r}} \men
    2^{\frac{1}{r}}\left( \max_{1 \men i \men k}\{1 - \al_i\} \right)^{\frac{2}{r}-1}
		\left(\sum_{i=1}^{k}{\al_i\norm{x_i}_p^r}\right)^{\frac{1}{r}}.
$$
The inequality (\ref{inequality3}) was deduced by S.A. Pichugov and V.I. Ivanov 
in (\cite{Pichugov2}, Assertion 1).

Using the Riesz-Thorin theorem for spaces with a mixed $L_p$-norm (\cite{Williams}, \S 14),  S.A.~Pichugov proved the following inequality 
(\cite{Pichugov1}, Theorem 1):
%The inequality (\ref{inequality1}) was proved  by Pichugov in \cite{Pichugov1},
\begin{gather} 
\nonumber
\left(\sum_{i=1}^{k}\sum_{j=1}^{l}{\al_i \beta_j \norm{(x_i - x_0) - (y_j - y_0)}_p^r}\right)^{\frac{1}{r}}  \\ \label{inequality4} \men
    2^{-\frac{1}{r'}} \left( \sum_{i_1=1,i_2=1}^{k}{\al_{i_1} \al_{i_2} \norm{x_{i_1} -x_{i_2}}_p^r} +
		\sum_{j_1=1,j_2=1}^{l}{\beta_{j_1} \beta_{j_2} \norm{y_{j_1} -y_{j_2}}_p^r}\right)^{\frac{1}{r}},
\end{gather}
	where  $\sum\limits_{i=1}^{k} {\al_i} =\sum\limits_{j=1}^{l} {\beta_j}=1, $ $ \al_i \bol 0 \;(i=1, \cdots, k) ,
	\; \beta_j \bol 0 \;(j=1, \cdots, l) ,$ $x_0 = \sum\limits_{i=1}^{k} {\al_i x_i},$ $y_0 = \sum\limits_{j=1}^{l} {\beta_j y_j}.$
	\\ \noindent Substituting $y_j$  for $0$ and $\beta_ j$ for $\frac{1}{l}$ in (\ref{inequality4}) we obtain the 
	inequality (\ref{inequality1}).
$\blacksquare$\\

{\bf Proof of Theorem \ref{th_Lp}.}\\
Consider  the case of $p \in (1, +\infty).$

For $L_p$ spaces and arbitrary set of vectors $A=\{x_0, x_1, \cdots, x_k\},$
 such that  $x_0 = \sum_{i=1}^{k} {\al_i x_i}, \; \sum_{i=1}^{k} {\al_i}=1, \; \al_i \bol 0 \;(i\in \ovl{1,k}), \; A\subset \B_1(0)$ we have
$$
    \left( \min_{i \in \ovl{1,k}}{\norm{x_0 -x_i}_p^r} \right)^\frac{1}{r} \men
     \left(\sum_{i=1}^{k}{\al_i \norm{x_i - x_0}_p^r}\right)^{\frac{1}{r}}.
$$

Using (\ref{inequality1}) and (\ref{inequality2}), since set of vectors $A$ was chosen arbitrarily,
we get  $\zeta_{L_p} \men 2^{ \left(\frac{1}{r} - \frac{1}{r'} \right)} = 2^{\left| \frac{1}{p} - \frac{1}{p'} \right|}$.

As it was shown in proof of Theorem \ref{UVOfinite} that $\zeta_{\ell_1^n}=\zeta_{\ell_{\infty}^n}=2\frac{n-1}{n}.$ 
Thus, $\zeta_{L_1}=\zeta_{L_\infty}=2.$ 
%Since Lebesgue integral is a limit of integrals of simple functions,
 %Lemma \ref{lemma_o_zamiuk} implies $\zeta_{L_p}= \lim\limits_{n \to \infty}{\zeta_{\ell^n_p}}.$ 
$\blacksquare$\\
\begin{remark}
       If $1 \men p \men 2,$ then,  using in the proof of Theorem \ref{th_Lp} inequality (\ref{inequality3}) instead of (\ref{inequality2}), we arrive at:
\begin{equation}\label{l_p_n_modul}
    \zeta_{\ell_p^n} \men \left(2\frac{n-1}{n}\right)^{\left| \frac{1}{p} - \frac{1}{p'} \right|}.
\end{equation}
\end{remark}

Still without any answer remains the following questions:
\begin{question}
Is  the  inequality (\ref{l_p_n_modul}) true if $p \in (2; \infty)$?
\end{question}
\begin{question}
Is the estimate in the  inequality (\ref{L_p_modul})  exact in case of $p \in (1; \infty), \, p \neq 2$?
\end{question}
%%%%%%%%%%%%%%%%%%%%%%%%%%%%%%%%%%%%%%%%%%%%%%%%%%%%%%%%
\section{Criterion of  a Hilbert space}
In order to prove Theorem \ref{UVOcriterion}  we need the following lemma, which follows directly from  the 
KKM theorem \cite{KKM_original}.
\begin{lemma} \label{lemma_triangle}
	Let $X$ be a Banach space.   Suppose the  triangle $a_1a_2a_3 \subset X$  satisfies the inequality  $ \diam{a_1a_2a_3} \men 2R$ and  is covered by balls $\B_R(a_i),\, i=1,2,3.$ Then these balls have a common point lying in the plane of the triangle.
\end{lemma}
%{\bf Proof.}\\
%Though there is an elementary (rather long) proof of this lemma, we provide the following short proof using the nerve of covering.
%	Let $L$ be a cross section of $X$ with a plane that contains the triangle $a_1a_2a_3.$ 
%	And $L\cap \B_R(a_i)= \B^2_R(a_i),$ where $\B^2_R(a_i) $ is a ball in Banach space $L$. 
%	The triangle $a_1a_2a_3$ is covered by balls $\B^2_R(a_i), i=1,2,3,$ 
%	so this covering is admissible. 
%	According to Lemma \ref{lemma_collapse} it is also contractible. 
%	Then the nerve of covering is contractible as well \cite{Nerv}. Since $ \diam{a_1a_2a_3} \men 2R,$ the balls intersect each other at every side of the triangle. Thus the nerve contains the boundary of the triangle and, since the triangle is contractible, it contains the triangle itself, so the balls intersection is nonempty.
%	$\blacksquare$

Taking into account Lemma \ref{lemma_triangle}, the proof of  Theorem \ref{UVOcriterion} is very similar to the one of Theorem 5 from   \cite{Garkavi}.

{\bf Proof of theorem \ref{UVOcriterion}.}\\
Using Theorem \ref{UVOfinite} and Corollary \ref{Hilbert} it suffices to prove that a Banach space $X,$ with $\dim X \bol 3$ and $\zeta_X=1$ is a Hilbert space. According to the well-known results obtained by Frechet and Blashke-Kakutani, it is enough to describe  only the case when $\dim X = 3.$ We need to show that if $\zeta_X=1,$ then for every 2-dimensional subspace there exists a unit-norm operator that projects $X$ onto this particular subspace. Let $0 \in L$ be an arbitrary 2-dimensional subspace in $X$, point $c$ is not contained in $L$.
We denote $\B_n^2(0)= L \cap \B_n(0)$ (it is a ball of radius $n \in \mathbb{N}$ in space $L$).
For every $n \in \mathbb{N}$ let us introduce the following notations:

\begin{eqnarray*}
    E_n = \left\{x\in L:\ \norm{c-x}\le n \right\},\\
    F_n = \left\{x\in L:\ \norm{c-x}= n \right\}.
\end{eqnarray*}

If $n$ is big enough, these sets are nonempty. 
Let $x_1,\, x_2, \, x_3$ be arbitrary points from $E_n.$ 
The CHD-constant of space $X$ equals 1, 
so the balls  $\B_n^2(x_i), i=1,2,3$ cover the triangle $x_1x_2x_3.$ 
According to Lemma \ref{lemma_triangle}, their intersection is not empty. 
According to Helly theorem, the set 
$$ S_n = \bigcap_{x \in E_n}{\B_n^2(x)}$$ 
is nonempty as well. 

Let us pick $a_n \in S_n,$ then by construction we have
\begin{equation}\label{crit1_1}
\norm{x-a_n} \men \norm{x -c}
\end{equation}
for every $x \in F_n.$
Let us show that
$$
\norm{x-a_n}\men \norm{x-c}
$$
 for every $x \in E_n.$
Suppose that for some $x \in E_n$
\begin{equation}\label{crit1_3}
 \norm{x-a_n} > \norm{x-c}.
\end{equation}

According to (\ref{crit1_1}) we may assume that $x \in E_n \setminus F_n.$
 Set $E_n$ is bounded and its boundary relatively to subspace $L$ coincides with $F_n$, thus there exists point $b \in F_n$, such that $x$ is contained in interval $(a_n,\;b).$ Then $a_n - x= \la(a_n - b), \; 0 <\la<1.$

Note that $ c- x= (c-a_n)+(a_n-x)=c-a_n + \la (a_n-b),$ 
then (\ref{crit1_3}) may be reformulated as $\norm{c-a_n + \la (a_n-b)} < \la \norm{a_n-b}.$

So,
\begin{eqnarray*}
\norm{c-b}=\norm{(c-a_n)+\la(a_n-b)+(1-\la)(a_n-b)}\men\\
\men \norm{(c-a_n)+\la(a_n-b)} + (1-\la)\norm{a_n-b} < \la \norm{a_n-b}+(1-\la)\norm{a_n-b}=\norm{a_n-b},
\end{eqnarray*}
and it contradicts (\ref{crit1_1}).

Consider the sequence $\{a_n\}.$ Note that $E_n\subset E_{n+1}$ and $\cup_{i=1}^{\infty}E_i = L.$ 
So, starting with a fixed natural $k$, 
the inclusion $0\in E_n, \; n \bol k$ becomes true, 
thus when $x=0$ inequality (\ref{crit1_1}) implies $\norm{a_n} \men \norm{c}, \; n \bol k,$
 i.e. the sequence $\{a_n\}$ is bounded. 
It means that sequence $\{a_{n}\}$ has a limit point $a$. 
Then every point $x \in L$ satisfies $\norm{x-a} \men \norm{x -c}.$ 
Let now represent every element $z\in X$ in the form
$$
    z = tc+x\;(x\in L, \, t \in \R).
$$

Operator $P(z)=P(tc+x)=ta+x$ projects $X$ onto $L.$

In addition:
$$
\norm{P(z)}=\norm{ta+x}=|t|\norm{a+\frac{x}{t}}\men|t|\norm{c+\frac{x}{t}}=\norm{tc +x}=\norm{z}.
$$

Hence, the $\norm{P}=1$ and taking into consideration the theorem of Blashke and Kakutani we come to a conclusion that $X$ is a Hilbert space.
$\blacksquare$\\

%%%%%%%%%%%%%%%%%%%%%%%%%%%%%%%%%%%%%%%%%%%%%%%%%%%%%%%%%%%%%%%%%%%%%%%%%%%%%%%%%%%%%%%%
{\bf Proof of Theorem \ref{ballcriterion}.}\\
It remains to check that in every Banach space $X$ that is not a Hilbert one, where $\dim{X}=3$, 
there exist a convex set and an admissible and not contractible covering.

To make the proof easier we first need  to prove a trivial statement from geometry.

Let hyperplane $H$ divide space $X$ in two half-spaces $H_+, H_-.$ 
Let $M$ be a bounded set in $H.$ We want to cover set $M$ with balls $B= \{\cup{\B_d(a_i)\mid i\in \ovl{1,n}, n \in \mathbb{N}}\}$ and call this covering  $(\e,r,H_+)$-good if $h^+(B,\; H_-) \men \e.$

\begin{lemma}\label{goodcover}
Let $X$ be a Banach space, $3 \men \dim X < +\infty$. Let hyperplane $H$ divide $X$ in two half-spaces: $H_+$ and  $H_-.$ Let $M$ be a bounded set in $H.$ Then for every $\e>0,\; d >0$ there exists an admissible set of balls ${\B_d(a_i), i\in \ovl{1,N}, N \in \mathbb{N}}$, such that set
 $B= \bigcup\limits_{i\in \ovl{1,N}}\B_d(a_i)$ 
may be regarded as $(\e,d,H_+)$-good covering of set $M$ and 
$\co\left(M\cup\{a_i\}\right) \subset B,\; i \in \ovl{1,N}.$
\end{lemma}
{\bf Proof.}\\
Let $\dim X = n.$
Without loss of generality we assume that $\e < d$ and $H$ is the supporting hyperplane for the ball $\B_d(0)$ 
and  $\B_d(0)\subset H_-.$ 
For any $r > 0$ and $a \in X$ we use $C_r(a)$ to denote a $(n-1)$-dimensional hypercube centered in a that lies 
in the hyperplane parallel to $H,$ where $r$ is the length of its edges.
Let $x \in H\cap\B_d(0).$ 
Then $h^+(\B_d(\e \frac{x}{\norm{x}}), H_-)\men  \e.$ Let $D= \B_d(\e \frac{x}{\norm{x}}) \cap H.$
 Note that $x$ is an inner point of set $D$ relatively to subspace $L$. 
In a finite dimensional linear space all norms are equivalent, so $C_r(x) \subset D$ for some  $r>0.$ 
As the ball $B_d(\e \frac{x}{\norm{x}})$ is centrally-symmetric,
it contains  affine hypercube $\co{\left(C_r(x)\cup C_r(\e \frac{x}{\norm{x}})\right)}$. 
Consider next an arbitrary bounded set $M \subset L.$ Since it is bounded, $M \subset C_R(b),$ where $b\in L,\; R>0.$ We suppose that $R= kr, \; k\in \mathbb{N}.$ Lets split hypercube $C_R(b)$ in hypercubes with edges of length $r$ and let $b_i, i\in \ovl{1,N}$ be the centers of these hypercubes. 
Hence, from the above, the balls $\B_d(b_i-(d-\e)\frac{x}{\norm{x}})$ give us the necessary covering.
$\blacksquare$

Let us consider an approach to construct  an admissible and not contractible covering of a convex set.

Let a Banach space $X$  be a non-Hilbert one and $\dim X = 3$. 
According to Theorem \ref{UVOcriterion}, $\zeta_{X} > 1,$  by Remark \ref{remark}, 
there exists set $A= \{a_1,\; a_2,\; a_3\} \subset \B_1(0)$ and point $b \in \co{A}$,
 such that $\rho(b,\; A)=1 +4\e > 1.$ According to Theorem \ref{UVOfinite}, $o \notin H.$
Consider the balls $\B_{1+\e}(a_i), \; i\in \ovl{1,3},$
 let $B_1 =\B_{1+\e}(a_1)\cup \B_{1+\e}(a_2)\cup\B_{1+\e}(a_3).$
 It is obvious that $b\notin B_1.$ Since all the edges of triangle $a_1a_2a_3$ lie in  $B_1,$ 
facets $0a_1a_2, \; 0a_1a_3,\; 0a_2a_3$ of tetrahedron $0a_1a_2a_3$ lie in $B_1.$ 
Let $H$ be a plane passing through points $a_1,\; a_2,\; a_3.$

Let $H$ divide space $X$ in two half-spaces: $H_+$ and $H_-.$ Let $0 \in H_+.$ 
According to Lemma \ref{goodcover} there exists an $(\e, 1+\e, H_+)$-good covering of triangle $a_1a_2a_3$ 
with an admissible set of balls 
that have centers lying in a set $C=\{c_i,\; i \in \ovl{1, N}\},$ $N\in \mathbb{N}.$ 
Let $B_2= \bigcup\limits_{i\in \ovl{1,N}}\B_{1+\e}(c_i).$
 Then set $B =B_1 \cup B_2$ cointains all the facets of tetrahedron $oa_1a_2a_3$ and does not contain interior of ball $\B_{\e}\left(b - 2\e\frac{b}{\norm{b}}\right),$ i.e. set $B$ is not contractible. However, $\co{\left(A\cup C\right)}\subset B_2 \subset B,$ i.e. union of balls $\B_{1+\e}(x), \; x \in A \cup C$ is an admissible covering for the set $\co{\left(A\cup C\right)}$ we were looking for.
$\blacksquare$

There remain still some open questions:

\begin{question} 
What is the minimal number of balls in an admissible and not contractible set of balls for a certain space $X$?  How to express this number in terms of space characteristics, such as its dimension, modulus of smoothness and modulus of convexity?
\end{question}

\begin{question}
How to estimate the minimal density of an admissible covering with balls for it to be contractible?
\end{question}

According to Lemma \ref{lemma_triangle}, it takes at least 4 balls to construct an admissible not contractible set of balls in an arbitrary Banach space. The following example describes the case with precisely 4 balls.
\begin{example}
    Let $X=l_1^3,$ $a_1 = \left(-\frac{2}{3},\frac{1}{3},\frac{1}{3}\right),\;
    a_2 = \left(\frac{1}{3},-\frac{2}{3},\frac{1}{3}\right),\;
    a_3 = \left(\frac{1}{3},\frac{1}{3},-\frac{2}{3}\right),$
     $a_4 = \left(-\frac{1}{6},-\frac{1}{6},-\frac{1}{6}\right).$
     Set of balls $\B_1(a_i), i =\ovl{1,4}$ is admissible, however, the complement of set
     $B=\bigcup\limits_{i=\ovl{1,4}}\B_1(a_i)$ has two connected components.
\end{example}

{\bf Proof}\\
1) Let us show that this set of balls is admissible. Every point $x$ from the tetrahedron $A=a_1a_2a_3a_4$ may be represented in form $x= \al_1 a_1+\cdots + \al_4 a_4,$ where $\al_1+ \cdots +\al_4=1, \; \al_i \bol 0, i\in\ovl{1,4}.$ Using equation $\al_4= 1 -\al_1-\al_2-\al_3,$ we are going to construct an inequation which would detect that point $x \in A$ is not contained in ball $\B_1(a_4):$

\begin{equation}\label{example}
 1 < \norm{x-a_4}=\left| \frac{-\al_1+\al_2+\al_3}{2} \right|+
 \left| \frac{\al_1-\al_2+\al_3}{2} \right|+ \left| \frac{\al_1+\al_2-\al_3}{2} \right|.
\end{equation}

Note that if every expression under the modulus  is positive,
 then the right side of (\ref{example}) equals $\frac{1}{2}(\al_1+\al_2+\al_3) \men \frac{1}{2}.$
 So, one of them has to be negative. 
Without lose of the generality, let $\alpha_1 \bol \alpha_2+\alpha_3.$ 
Then the other two expressions are positive and inequation (\ref{example}) can be rewritten:
$
    \frac{3\al_1-\al_2-\al_3}{2}>1.
$
Then $\al_1 > \frac{2}{3}+ \frac{1}{3}(\al_2+\al_3).$

Using this relation we arrive at:
$$
    1 -\al_4 =\al_1+\al_2+\al_3\bol \frac{2}{3}+\frac{4}{3}(\al_2+\al_3).
$$

Thus,
$
    \frac{1}{4} - \frac{3}{4}\al_4 \bol \al_2+\al_3.
$

We use the last inequality to estimate the distance between $x$ and the vertex $a_1:$
\begin{gather*}
    \norm{x-a_1}= \norm{\al_2(a_2-a_1)+\al_3(a_3-a_1)+\al_4(a_4-a_1)}\\
    \men \al_2\norm{a_2-a_1}+\al_3\norm{a_3-a_1}+\al_4\norm{a_4-a_1} \\
    =2(\al_2+\al_3) +\frac{3}{2}\al_4\men 2(\frac{1}{4} - \frac{3}{4}\al_4)+\frac{3}{2}\al_4 =\frac{1}{2}.
\end{gather*}

So, we come to the conclusion that the  set of balls is admissible. 	

2)Let $b_1= \left(\frac{1}{3},\frac{1}{12}, \frac{1}{12}\right), \;
b_2= \left(\frac{1}{12},\frac{1}{3}, \frac{1}{12}\right),\;
b_3= \left(\frac{1}{12},\frac{1}{12}, \frac{1}{3}\right),\;
b_4= \left(\frac{1}{3},\frac{1}{3}, \frac{1}{3}\right),$
tetrahedron $\Delta = b_1b_2b_3b_4.$

It is easy enough to show that $\partial \Delta  \subset B,$ but $\vn{ \Delta} \cap B = \emptyset.$
$\blacksquare$\\
%%%%%%%%%%%%%%%%%%%%%%%%%%%%%%%%%%%%%%%%%%
\section{About contractibility of a proximally smooth sets} 
Clark, Stern and Wolenski \cite{Clark} introduced and studied the {\it proximally smooth sets} in a Hilbert space 
$H.$ 
A set $A \subset X$ is said to be proximally smooth 	with constant $R$ if the distance function
 $x \rightarrow \rho(x, A)$ is continuously differentiable on  set 
$U(R,A)=\left\{ x \in X: 0 < \rho(x, A) < R \right\}.$
Properties of proximally smooth sets in a Banach  space and relations between such sets and akin classes of set, including uniformly prox-regular sets, were investigated in \cite{Clark}-\cite{t1}. We study the sufficient condition of the contractibility for a proximal smooth sets.
G.E. Ivanov showed that if $A \subset H$ is proximally smooth (weakly convex in his terminology) with constant $R$ 
and $A \subset B_r(o)$ with $r < R,$ then $A$ is contractible.
The following theorem is a generalization of this result.
%Согласно работе \cite{Clark} множество $A \subset X$ называется {\it
%проксимально гладким} с константой $R$, если функция расстояния
%$x \rightarrow \rho(x, A)$  непрерывно дифференцируема на
%множестве $U(R,A)=\left\{ x \in X: 0 < \rho(x, A) < R \right\}$.

%В работе \cite{IvBal} показано, что в равномерно выпуклом и равномерно
%гладком банаховом пространстве $X$ метрическая проекция на
%замкнутое проксимально гладкое с константой $R$ множество
%$A\subset X$ непрерывна на множестве $U(R,A)$. Отсюда и изы
%теоремы 1 получаем следующий результат.

\begin{fact}\label{th_prox}
Let $X$ be a uniformly convex and uniformly smooth Banach space.
Let $A$ be a closed and proximally smooth with constant $R$ subset of a ball with radius  $r<\frac{R}{\zeta_X}.$ Then $A$ is contractible.
 %Пусть замкнутое множество $A$  из
%равномерно выпуклого и равномерно гладкого банахова пространства
%$X$ является проксимально гладким с константой $R$ и  содержится в
%шаре радиуса $r<\frac{R}{\zeta_X}$. Тогда $A$ стягиваемо.
\end{fact}
{\bf Proof.}\\
Note that set $\co A$  is contractible, so a continuous function 
$F:[0,1]\times\co A\to\co A$ such that  
$F(0,x)=x$, $F(1,x)=x_{0}$ for all $x\in \co A$ and some $x_0 \in A$ exist.
Due to the CHD-constant definition and inequality $r < \frac{R}{\zeta_X}$ the set $\co A$ belongs to  the $R$-neighborhood of the set $A.$
On the other hand,	$A$ is proximally smooth and in accordance to paper \cite{IvBal} metric projection mapping 
$\pi:\co A\to A$
is single valued and continuous.
Finally, we define the mapping $\tilde F:[0,1]\times A\to A$  as follows  $\tilde
F(t,x)=\pi(F(t,x))$ for all $t\in[0,1]$, $x\in A.$ 
The mapping $F$ contracts  set $A$ to point $x_0$.
%стягиванием множества $A$.
%Заметим, что поскольку множество $\co A$ выпукло и ограничено, то
%оно стягиваемо, то есть, существует точка $x_{0}\in\co A$ и
%непрерывная функция $F:[0,1]\times\co A\to\co A$ такие, что
%$F(0,x)=x$, $F(1,x)=x_{0}$ для любого $x\in \co A$. Из определения
%УВО-модуля следует, что множество $\co A$ содержится в  $R$-окрестности множества $A$. С другой стороны, $A$ является
%проксимально гладким с константой $R$ множеством, а значит,
%согласно работе \cite{IvBal}, отображение метрического проектирования
%$\pi:\co A\to A$ однозначно и непрерывно. Поэтому отображение
%$\tilde F:[0,1]\times A\to A$, заданное формулой $\tilde
%F(t,x)=\pi(F(t,x))$ при всех $t\in[0,1]$, $x\in A$, является
%стягиванием множества $A$.
$\blacksquare$

%\bibliographystyle{model1-num-names}
%\bibliographystyle{ieeetr}
%\bibliography{uvolit}

\end{document}